\documentclass{amsart}

\usepackage{amssymb}
\newtheorem{theorem}{Theorem}

\title{Cancellation does not imply stable rank one}

\author{Andrew S. Toms}

\thanks{This work was supported by an NSERC Postdoctoral Fellowship}

\begin{document}
\maketitle

\begin{abstract}

An unital $C^*$-algebra $A$ is said to have cancellation of projections if the semigroup 
$D(A)$ of Murray-von Neumann equivalence classes of projections in matrices over $A$ is
cancellative. It has long been known that stable rank
one implies cancellation for any $A$, and some partial converses have been 
established.  In the sequel it is proved that cancellation does not imply stable rank one
for simple, stably finite $C^*$-algebras.

\end{abstract}

\section{Introduction}

Rieffel introduced the notion of stable rank for $C^*$-algebras in his 1983 paper \cite{R}:  a 
unital $C^*$-algebra $A$ is said to have stable rank $n$ ($\mathrm{sr}(A) = n$) if $n$ is the least
natural number such that the set
\[
\mathrm{Lg}_n(A) \stackrel{\mathrm{def}}{=} \left\{ (a_1,\ldots,a_n)\in A^n| \exists b_i \in A_i, 
\ 1 \leq i \leq n \ : \ \sum_{i=1}^n b_i a_i = 1 \right\}
\]
is dense in $A^n$.  If no such $n$ exists, then one says that the stable rank of $A$  
is infinite.  In the case of a commutative $C^*$-algebra, the stable rank is proportional
to the covering dimension of the spectrum;  stable rank may be viewed as a kind of 
non-commutative dimension. 

Given an unital $C^*$-algebra $A$, let $D(A)$ be the Abelian semigroup obtained by 
endowing the set of Murray-von Neumann equivalence classes of projections in matrix 
algebras over $A$ with the addition operation coming from direct sums.  The algebra 
$A$ is said to have cancellation of projections if $x+y=x+z$ implies that $y=z$ for any $x,y,z \in D(A)$.
Shortly after the appearance of Rieffel's paper, Blackadar showed that stable rank one
 implies cancellation of projections (\cite{B}).  He also established a partial
converse:  if a $C^*$-algebra of real rank zero has cancellation of projections, then
it has stable rank one.  The relationship between cancellation and stable rank for 
general simple, stably finite $C^*$-algebras, however, remained unclear. 
The lack of examples of simple, stably finite $C^*$-algebras with non-minimal stable rank was a serious obstacle.  
Villadsen provided the first such examples in \cite{V}, but determining 
whether his examples had cancellation of projections was all but impossible due to their extremely complicated 
$\mathrm{K}$-theory.

Recently, the author has been able to apply Villadsen's techniques to construct simple, stably
finite $C^*$-algebras with non-minimal stable rank and cyclic $\mathrm{K}_0$-groups.  These
algebras constitute the first simple, nuclear and stably finite counterexample to Elliott's 
classification conjecture for nuclear $C^*$-algebras (\cite{E}, \cite{T}).  In the sequel
we study one such algebra in order to prove our main result.
\begin{theorem}
There is a simple, separable, nuclear, and stably finite $C^*$-algebra with non-minimal
stable rank which nevertheless has cancellation of projections.
\end{theorem}

\noindent
Thus, Blackadar's partial converse cannot be extended to cover general simple, stably finite
$C^*$-algebras. 

\vspace{2mm}
\noindent 
{\it Acknowledgements.}  The author would like to thank the referee for several comments which
improved the exposition of the sequel.  

\section{The proof of the main result}

\begin{proof}
We proceed by a close analysis of the structure of
the simple, separable, and stably finite $C^*$-algebra $B_2$ of 
\cite{T}, which has non-minimal stable rank.  We will prove that
$B_2$ nevertheless has cancellation of projections. 

Let $C$, $D$ be $C^{*}$-algebras and let $\phi_{0}$, $\phi_{1}$ be
$*$-homomorphisms from $C$ to $D$.  The generalised mapping torus of
$C$ and $D$ with respect to $\phi_{0}$ and $\phi_{1}$ is 
\begin{displaymath}
A:=\{(c,d)|d \in C([0,1];D), \ c \in C, \ d(0)=\phi_{0}(c), \ d(1)=\phi_{1}(c)\}
\end{displaymath}
We will denote $A$ by $A(C,D,\phi_{0}, \phi_{1})$ for clarity when necessary.   
Let $\mathcal{U}(A)$ denote the unitary group of an unital $C^*$-algebra $A$.

The algebra $B_2$ of \cite{T} is constructed as the limit of an inductive sequence
$(A_i,\theta_i)$ of generalised mapping tori $A_i = A(C_i,D_i,\phi_i^0,\phi_i^1)$ and
unital $*$-homomorphisms $\theta_i:A_i \to A_{i+1}$ where, for each $i \in \mathbb{N}$, 
\[
C_i \stackrel{\mathrm{def}}{=} p_i (\mathrm{C}(X_i) \otimes \mathcal{K}) p_i
\]
and 
\[
D_i \stackrel{\mathrm{def}}{=} \mathrm{M}_{k_i} \otimes C_i
\]
for some connected compact Hausdorff space $X_i$, projection $p_i \in
 \mathrm{C}(X_i) \otimes \mathcal{K}$ and natural number $k_i$.  The maps
$\phi_i^0$ and $\phi_i^1$ are unital.  The spaces $X_i$, $i \in \mathbb{N}$, have the property that 
\[
\mathrm{dim}(p_i) = \frac{\mathrm{dim}(X_i)}{2},
\]
and the maps $\phi_i^0$ and $\phi_i^1$ are chosen to ensure that
\[
\left(\mathrm{K}_0A_i,\mathrm{K}_0A_i^+,[1_{A_i}]\right) = (\mathbb{Z},\mathbb{Z}^+,1),
\]
where $1_{A_i} \in A_i$ is the unit; $A_i$ is projectionless but for zero and $1_{A_i}$.

To prove Theorem 1 it will suffice to prove that $A_i$ has
cancellation of projections for every $i \in \mathbb{N}$.
Let $p,q \in \mathrm{M}_n(A_i)$ be projections having the same
$\mathrm{K}_0$-class.  We must show that $p$ and $q$ are Murray-von Neumann
equivalent.    Since
$\mathrm{K}_0(A_i) = \mathbb{Z}[1_{A_i}]$, we may assume that
$p$ is a multiple of the unit of $A_i$, say $p = l 1_{A_i}$. 
$\mathrm{M}_n(A_i)$ can be viewed as an algebra of functions from $[0,1] \times X_i$
into matrices.  Given $f \in \mathrm{M}_n(A_i)$, we let $f(t)$, $t \in [0,1]$, denote the restriction
of $f$ to $\{t\} \times X_i \subseteq [0,1] \times X_i$.  Both $f(0)$ and $f(1)$ are
images of a single element in $\mathrm{M}_n(C_i)$, which we denote by $f(\infty)$. 
If two vector bundles over a compact, connected
CW-complex $X$ of covering dimension $m$ with the same $\mathrm{K}^0$-class
have fibre dimension at least $m/2$, then the bundles are isomorphic (cf. Theorem 1.5, Chapter 8,
\cite{H}).  In the language of $C^*$-algebras, the projections in $\mathrm{M}_k \otimes \mathrm{C}(X)$,
some $k \in \mathbb{N}$, corresponding to these vector bundles are Murray-von Neumann
equivalent.  
Since $p(\infty)$ and $q(\infty)$ can be viewed as vector bundles over $X_i$ having the same
$\mathrm{K}^0$-class, and since they must both have fibre dimension at least
$\mathrm{dim}(X_i)/2$ by the construction of $A_i$, they are Murray-von Neumann equivalent, as
are their images under $\phi_i^0$ and $\phi_i^1$.  Note that if one considers
$\mathrm{M}_n(A_i)$ as an unital sub-$C^*$-algebra of \mbox{$C_i \otimes \mathrm{M}_{n k_i} \otimes 
\mathrm{C}([0,1])$}, then fibre dimension considerations show $q$ and $p$ to be Murray-von Neumann
equivalent inside $C_i \otimes \mathrm{M}_{n k_i} \otimes \mathrm{C}([0,1])$.
This does not, however, prove that $q$ and $p$ are
Murray-von Neumann equivalent inside $\mathrm{M}_n(A_i)$.    

We may assume without loss of generality
that $l 1_{A_i}$ and $q$ are constant over some small interval $[\frac{1}{2} - \epsilon, 
\frac{1}{2} + \epsilon]$ in the interval factor of the spectrum of $\mathrm{M}_n(A_i)$, 
since small perturbations do not disturb the 
Murray-von Neumann equivalence class.  Consider $l 1_{A_i}$ and $q$
as vector bundles over $[0,1] \times X_i$.  Define
\[
q_0 := q|_{[0,\frac{1}{2}-\epsilon] \times X_i}, \ \
q_1 := q|_{[\frac{1}{2}+\epsilon,1] \times X_i}
\]
and
\[
1_{A_i,0} := 1_{A_i}|_{[0,\frac{1}{2}-\epsilon] \times X_i}, \ \
1_{A_i,1} := 1_{A_i}|_{[\frac{1}{2}+\epsilon,1] \times X_i}.
\]

Corollary 4.4, Chapter 3, \cite{H}, states:

\vspace{1mm}
\begin{enumerate} 
\item[] Let $\gamma$ be a vector bundle over
$X \times [0,1]$, $X$ paracompact, and $\omega$ a vector
bundle over $X$ such that $\gamma|_{X \times \{0\}} \cong \omega$.
Then, $\gamma$ is isomorphic to the induced bundle $\pi^{*}(\omega)$, 
where $\pi:X \times [0,1] \to X \times \{0\}$ is given by $\pi(x,t) = (x,0)$. 
\end{enumerate}

\vspace{1mm}
\noindent
Define maps
\[
\pi_0:[0,1/2-\epsilon] \times X_i \to \{0\} \times X_i, \ \ 
\pi_1:[1/2+\epsilon,1] \times X_i \to \{1\} \times X_i
\]
by
\[
\pi_0(t,x) = (0,x), \ \ \pi_1(t,x) = (1,x).
\]
We have $l1_{A_i}(j) \cong q(j)$ for $j \in \{0,1\}$.  Moreover,
$l 1_{A_i,j} \cong \pi_j^*(l 1_{A_i}(j))$ by construction.  
We may thus apply Corollary 4.4, Chapter 3, \cite{H},
with $\gamma = q_j$, $\omega = l 1_{A_i}(j)$, and $\pi = \pi_j$ 
to conclude that $l 1_{A_i,j} \cong q_j$.  In other words, there is a continuous path of partial isometries $v(t)$, $t \in
[0,\frac{1}{2}-\epsilon] \cup [\frac{1}{2}+\epsilon,1]$, such that $v(t)^*v(t) = l1_{A_i}(t)$,
 $v(t)v(t)^* = q(t)$, and, for each $j \in \{0,1\}$, the partial isometry $v(j)$
is the image under $\phi_i^j \otimes \mathrm{id}_{\mathrm{M}_n}$ of a single partial isometry
$v \in \mathrm{M}_n(C_i)$ such that $v^*v = l1_{C_i}$ and $vv^* = q(\infty)$.
This last property ensures that if we can find a continuous extension of $v(t)$ to a partial 
isometry defined on $[0,1]$, then our proof is complete --- $v(t)$ will lie in $\mathrm{M}_n(A_i)$.

From \cite{N} we have the formula
\[
\mathrm{sr}(p(\mathrm{C}(X) \otimes \mathcal{K})p) = \left\lceil 
\frac{\lfloor \mathrm{dim}(X)/2 \rfloor}{\mathrm{rank}(p)} \right\rceil + 1,
\]
where $\mathcal{K}$ denotes the compact operators on a separable Hilbert space, $X$ is 
a compact connected Hausdorff space, and $p$ is a projection in $\mathrm{C}(X) \otimes \mathcal{K}$.
Straightforward calculation then shows that $\mathrm{sr}(C_i)=2$, $\forall i \in \mathbb{N}$. 
For an unital $C^*$-algebra $A$, let $\mathcal{U}(A)$ denote the unitary group of $A$, and 
let $\mathcal{U}(A)_0$ denote the connected component of $\mathcal{U}(A)$ containing the identity.  Theorem
10.12 of \cite{R} states that one has an isomorphism
\[
\frac{\mathcal{U}(\mathrm{M}_r(A))}{\mathcal{U}(\mathrm{M}_r(A))_0} \to \mathrm{K}_1(A)
\]
whenever $r \geq \mathrm{sr}(A)+2$.  In the construction of $A_i$, the parameter $k_i$ in 
the definition $D_i := \mathrm{M}_{k_i}(C_i)$ is chosen to be much larger than $\mathrm{sr}(C_i)$.
Furthermore, one has (again, by construction) that $\mathrm{K}_1(C_i) = 0$, $\forall i \in \mathbb{N}$.
Thus, $\mathcal{U}(\mathrm{M}_l(D_i))$ is connected for every $l \in \mathbb{N}$. 

We may view $u:=v(1/2+\epsilon)^*v(1/2-\epsilon)$ as a unitary element in $\mathrm{M}_l(D_i)$.  By the
discussion above, there is a path of unitary elements $u(t)$, $t \in [1/2-\epsilon,1/2+\epsilon]$, 
inside $\mathrm{M}_l(D_i)$ such that $u(1/2+\epsilon) = l1_{A_i}$ and $u(1/2-\epsilon) =u$.
For $t \in [1/2-\epsilon,1/2+\epsilon]$, define $\tilde{v}(t) = v(1/2+\epsilon)u(t)$.  Clearly, $\tilde{v}(t)$
is a partial isometry in $\mathrm{M}_n(D_i)$ for each $t$ in its domain.  One has 
\[
\tilde{v}(1/2+\epsilon) = v(1/2+\epsilon)
\]
and
\[
\tilde{v}(1/2-\epsilon) = v(1/2+\epsilon)v(1/2+\epsilon)^*v(1/2-\epsilon) = q(1/2-\epsilon)v(1/2-\epsilon) = v(1/2-\epsilon).
\]
Then
\[
v(t) := \left\lbrace \begin{array}{ll} v(t), \  t \in [0,\frac{1}{2}-\epsilon] \cup [\frac{1}{2}+\epsilon,1]\\
\tilde{v}(t), \ t \in (1/2-\epsilon,1/2+\epsilon) \end{array} \right.
\]
defines a partial isometry in $\mathrm{M}_n(A_i)$ such that $v(t)^*v(t) = l1_{A_i}(t)$ and $v(t)v(t)^* = q(t)$,
$\forall t \in [0,1]$. $q$ and $l1_{A_i}$ are thus Murray-von Neumann equivalent, as desired.

\end{proof}

\vspace{5mm}

\noindent
\emph{Andrew S. Toms \newline
Department of Mathematics and Statistics \newline
\hspace*{2mm} University of New Brunswick \newline
\hspace*{2mm} Fredericton, New Brunswick, E3B 5A3 \newline
Canada \newline}

\noindent
atoms@unb.ca

\end{document}